\documentclass[12pt]{article}
\usepackage{amsfonts,amssymb,amsthm,amsmath}
\theoremstyle{theorem}
\newtheorem*{Th}{Theorem}
\theoremstyle{lemma}
\newtheorem{Lm}{Lemma}
\newcommand{\probel}{\mbox{ }}
\newcommand{\ir}{[0,\vartheta]\times \mathbb{R}^n}

\begin{document}
\title{A Transformation of the Control under Uncertainty Problems}
\author{Yurii Averboukh}
\date{}
\maketitle
\begin{abstract}
 The game theoretical approach problem is considered. If goal set is controllability set of auxiliary control system then the original problem can be transformed into the problem of approach ``at the moment''.
\end{abstract}

\section{Introduction}
The paper is devoted to the theory of differential games. Within the framework of this theory the control processes  with conflict or uncertainty are studied. The investigation of such problem started with the book of R. Isaacs \cite{isaacs}. Basic approaches to the mathematical theory of conflict controlled processes were obtained by L.S.~Pontryagin, N.N.~Krasovski and B.N.~Pshenichnyi. The construction of the strong mathematical theory of differential games is connected with the research of N.N.~Krasovskii and A.I.~Subbotin \cite{NN_PDI_en}.
They have proved the alternative theorem. This theorem fully describes the structure of differential game.
The differential game can be reduced to the series of ordinary control problems due to A.G.~Chentsov Programmed Iteration Method (see \cite{Subb_Ch}--\cite{Chentsov_76}). This reduction widely uses in this paper.

The important class of differential games is the class of problem of approach ``at the moment''. In particular these problems can be studied with the help of theory of minimax (or viscosity) solutions constructed by A.I.~Subbotin (see \cite{Subb_book}). In this paper the method of reduction of wide class of differential games to the problem of approach ``at the moment'' is introduced. The reduction is realized by the substitution of control spaces and dynamic function. Initial and transformed games are equivalent in the sense of Programmed Iteration Method.

\section{Definitions and Preliminaries}
Consider a differential games in which the motion of controlled system is governed by ordinary differential equation on the segment $t\in [0,\vartheta]$
\begin{equation}\label{sys} \dot{x}=f(x,u,v),\probel x\in\mathbb{R},\probel u\in
P, v\in Q.
\end{equation}
Here $u\in P$ and $v\in Q$ are the controls of first player and second player respectively. The first player tries to bring the system onto the set $M$, $M\subset\ir$. The aim of the second player is opposite. The problem of first player is often called $M$-approach problem.

The sets $M$, $P$, $Q$ and function $f$ satisfy the following assumptions.
\begin{enumerate}
  \item The target set $M$ is closed.
  \item $P\subset \mathbb{R}^p$, $Q\subset\mathbb{R}^q$ for appropriate natural numbers $p$ and $q$.
  \item The sets $P$ and $Q$ are compact.
  \item $f(\cdot,\cdot,\cdot)$ is continuous.
  \item $f(\cdot,\cdot,\cdot)$ is locally lipschitzian with respect to the phase variable.
  \item $f(\cdot,\cdot,\cdot)$ satisfies sublinear growth condition.
\end{enumerate}
The differential game is considered in the class of contrstrategies of the first player and positional strategies of the second player.

According to the formalization of differential game suggested by N.N.~Krasovskii and A.I.~Subbotin any function $U:\ir\times Q\rightarrow P$ measurable with respect to the 3-d argument is called a contrstrategy of the first player, function $V:\ir\rightarrow Q$ is called a strategy of the second player. Let us define step-by-step and constructive motions. These definitions follow the book \cite{NN_PDI_en}. Let $U$ be a contrstrategy of the first player, $v(\cdot)$ be measurable control of the second player, $(t_*,x_*)$ be a position, $\Delta=\{\tau_k\}_{k=0}^{r}$ be a partition of the segment $[t_*,\vartheta]$. The function which satisfy the conditions
$$x[t]=x[\tau_{k-1}]+\int_{\tau_{k-1}}^{t}f(x[\theta],U(\tau_{k-1},x[\tau_{k-1}],v(\theta)),v(\theta)
)d\theta,\probel \probel t\in [\tau_{k-1},\tau_k] \forall k\in \overline{1,r},\probel x[\tau_0]\triangleq x_*
$$ is called step-by-step motion. Obviously, such function exists and it is unique. The limits of step-by-step motions as fineness of partition goes to 0 are called constructive motions in sense of N.N.~Krasovskii and A.I.~Subbotin. It is supposed that first player tries to bring all of constructive motions on the target set. The control of the second player is supposed to be formed by the rule $$v(t)=V(\xi_{i-1},x[\xi_{i-1}])\ \ \forall t\in [\xi_{i-1},\xi_i). $$ Here $\{\xi_i\}_{i=0}^m$ is a partition of the segment $[t_*,\vartheta]$.

By the alternative theorem of N.N.~Krasovskii and A.I.~Subbotion \cite{NN_PDI_en} the solution of approach problem is completely determined by the set of successful solvability. First player can bring the motion onto the target set iff it begins on the position from the solvability set. Denote this set by $\mathfrak{W}$. This set is maximal $u$-stable bridge. Moreover the solving contrstrategy is defined by the extremal shift rule.

The set ${W}$ is called $u$-stable bridge if for all $v\in Q$ and for all $(t_*,x_*)\in W$ there exists solution of differential inclusion $$\dot{y}(t)\in{\rm co}\{f(t,x,u,v):u\in P\} $$ $y(\cdot)$ and moment $\xi\in [t_*,\vartheta]$ such that $y(\xi)\in M[\xi]$ and for all $\in [t_*,\xi]$ the following inclusion holds: $y(t)\in W[t]$.

If the Isaacs condition $$\forall s,x\in \mathbb{R}^n\probel \min_{u\in P}\max_{v\in Q}\langle s,f(x,u,v)\rangle= \max_{v\in Q}\min_{u\in
P}\langle s,f(x,u,v)\rangle$$ is fulfilled, one can consider the differential game in the class of positional strategies (see \cite{NN_PDI_en}).

Ordinary control systems are utilized in this paper also. Consider a system
\begin{equation}\label{g_control}
\dot{x}=h(x,b),\probel b\in\Lambda.
\end{equation}
Denote by $S_{h,b}^\tau$  the flow for time $\tau$ generated by the constant control $b\in \Lambda$. We suppose that $\tau\in\mathbb{R}$. Let $(t_*,x_*)\in\ir$, $b(\cdot):[0,\vartheta]\rightarrow \Lambda$ be a measurable function, the solution of the equation $$x(t)=x_*+\int_{t_*}^th(x(\xi),b(\xi))d\xi $$ is called the motion of system  (\ref{g_control}) generated by the control $b(\cdot)$ and is denoted by  $x_h(t,t_*,x_*,b(\cdot))$. We assume that $t\in [0,\vartheta]$.

Let us consider the case when $b(\cdot)$ is piecewise constant function. Suppose that $t\geq t_*$. There exist the collection of numbers $\tau_1,\ldots,\tau_k\in [0,\vartheta]$ and the collection of controls $b_1,\ldots,b_k\in \Lambda$ such that $t=t_*+\tau_1+\ldots+\tau_k$ and $b(\xi)=b_i$  for $\xi\in [t_*+\tau_1+\ldots+\tau_{i-1},t_*+\tau_1+\ldots+\tau_{i-1}+\tau_i)$. In this case the following representation is fulfilled:
$$x_h(t,t_*,x_*,b(\cdot))=S_{h,b_k}^{\tau_k}\circ\ldots\circ S_{h,b_1}^{\tau_1}(x_*). $$ One can write an analogous representation for the case $t\leq t_*$. In this case $\tau_i\leq 0$.

The slide controls (or measure controls) are very useful in the control theory. Consider the Borel $\sigma$-algebra of subsets of  $[0,\vartheta]\times\Lambda$. Denote the set of measures defined on this $\sigma$-algebra by $\mathcal{R}_\Lambda$. In the control theory the elements of $\mathcal{R}_\Lambda$ are called slide controls. Let $\mu\in\mathcal{R}_\Lambda$, $(t_*,x_*)\in\ir$. The motion generated by slide control $\mu$ is the solution of equations $$
x(t)=x_*+\int_{[t_*,t]\times\Lambda}h(x(\xi),b)\mu(d(\xi,b))\mbox{ for }t\geq t_*,
$$
\begin{equation}\label{formula_motion_minus}
  x(t)=x_*-\int_{[t,t_*]\times\Lambda}h(x(\xi),b)\mu(d(\xi,b))\mbox{ for }t\leq t_*.
\end{equation}
note this motion by $\varphi_h(\cdot,t_*,x_*,\mu)$. If $\mu\in\mathcal{R}_\Lambda$ then there exist a sequence of piecewise constant controls $\{b^k(\cdot)\}_{k=1}^\infty$ such that $$x_h(\cdot,t_*,x_*,b^k(\cdot))\rightrightarrows \varphi_h(\cdot,t_*,x_*,\mu),\probel k\rightarrow\infty.$$ (See for details \cite{NN_Cyber}.) The set of ordinary controls can be embedded into the set of slide controls. Namely, let $b(\cdot):[0,\vartheta]\rightarrow \Lambda$ be a measurable control, there exists a measure $\mu_{b(\cdot)}\in\mathcal{R}_\Lambda$ such that $$\int_{[0,\vartheta]\times \Lambda}\psi(t,b)\mu_{b(\cdot)}(d(t,\mu))=\int_{0}^\vartheta\psi(t,b(t))dt $$ for all $\psi\in C([0,\vartheta]\times\Lambda)$. Note that $$x_h(\cdot,\cdot,\cdot,b(\cdot))=\varphi_h(\cdot,\cdot,\cdot,\mu_{b(\cdot)}).$$

The problem of differential game can be reduced into the series of control problem due to the Programmed Iteration Method suggested by A.G.~Chentsov \cite{Subb_Ch}. Let us define the program absorption operator. This operator is defined on the family of closed subset of $\ir$. Let us consider the controlled system (\ref{sys}). For any $v\in Q$ define the ordinary control system by the rule $$f_v(x,u)\triangleq f(x,u,v).$$ Let the the program absorption operator $A$ be given by
\begin{multline*}
E\mapsto A_f(E)\triangleq\{(t_*,x_*)\in E:\forall v\in Q\probel\exists \mu\in \mathcal{R}_P\probel\exists \xi\in [t_*,\vartheta]:\\ (\varphi_{f_v}(\xi,t_*,x_*,\mu)\in M[\xi])\ \ \&\ \ (\varphi_{f_v}(t,t_*,x_*,\mu)\in E[t]\probel\forall t\in [t_*,\xi])\}.
\end{multline*}
Here $E\subset\ir $, $E$ is closed, $$E[t]\triangleq\{x:(t,x)\in E\}. $$

Let us consider the sequence $$W_0\triangleq\ir,\probel W_k=A_{f}(W_{k-1}),\probel \forall k\in\mathbb{N}. $$  A.G. Chentsov established that the set of successful solvability of the first player problem can be represented in the form
\begin{equation}\label{Chentsov_formula}
\mathfrak{W}=\bigcap_{k=0}^\infty W_k.
\end{equation}

Note that the u-stability condition can be written in the terms of program absorption operator: i.e. the set $W$ is $u$-stable bridge iff $A(W)=W$.

\section{Main result}

Consider the differential game with controlled system (\ref{sys}) and target set $M\subset\ir$. Denote $F\triangleq M[\vartheta]$. Suppose that $M$ is controllability set of control system
$g(x,\omega)$, $\omega\in\Omega$, and the target set $M^*\triangleq \{\vartheta\}\times F$:
  \begin{equation}\label{M_control}
  M=\{(t,x)\in \ir:\exists x_*\in F \ \ \exists \mu\in\mathcal{R}_\Omega:\probel x=\varphi_g(t,\vartheta,x_*,\mu)\}.
\end{equation}
Here $\varphi_g(t,\vartheta,x_*,\mu)$ is defined by (\ref{formula_motion_minus}) It is assumed that $\Omega$ is compact set in finitely dimensional euclidean space.

Below we introduce transformed differential game ``at the moment''. In this game the resources of the first player are expanded by adding the control parameters $\nu$ and $\omega$. Consider on the segment $[0,\vartheta]$ the controlled system \begin{equation}\label{eq_sys}
    \dot{x}=f^*(x,\nu,u,\omega,v), \probel  x\in \mathbb{R}^n,\probel \nu\in \{0,1\}, \probel u\in P,\probel \omega \in \Omega,\probel v\in Q.
\end{equation}  Here
\begin{equation}\label{f_star}
f^*(x,\nu,u,\omega,v)=\nu\cdot f(x,u,v)+(1-\nu)\cdot g(x,\omega)=
\left\{
\begin{array}{lc}
  f(x,u,v), & \nu=0, \\
  g(x,\omega), & \nu=1.
\end{array}
\right.
\end{equation}

In system (\ref{eq_sys}) the variables $\nu$, $u$ and $\omega$ are controls of the first player, variable $v$ is control of the second player. Consider the $M^*$-approach problem for the system (\ref{eq_sys}). Analogous methods of control problem transformation were used in the case when target set is cylinder (see \cite{Agrachev}, \cite{Mitchel}).

Let us introduce the following notation. $${P}^*=\{0,1\}\times P\times \Omega.$$ If $h=f_v$, $u\in P$, then we denote $S^\tau_{h,u}$ by $\mathcal{F}_{u,v}^\tau  $; if $h={f}^*_{v}$, ${u}^*=(\nu,u,\omega)\in {P}^*$, then we denote $S^\tau_{h, u^*}$ by ${\mathcal{F}}_{{u}^*,v}^{*\probel\tau}$. Further denote $\mathcal{G}_\omega^\tau\triangleq S_{g,\omega}^\tau$.  Let $(t_*,x_*)\in\ir$, $t\in [0,\vartheta]$. If $h=f_v$, $\mu\in \mathcal{R}_P$, then we denote $\varphi_h(t,t_*,x_*,\mu)$ by $\phi(t,t_*,x_*,\mu,v)$. Analogously in the case $h=f^*_v$, $\mu\in\mathcal{R}_{P^*}$ $\varphi(t,t_*,x_*,\mu)$ is denoted by $\phi^*(t,t_*,x_*,\mu,v)$. Let $A$ be the program absorption operator for the initial problem, $A^*$ be  the program absorption operator for the transformed problem. Denote $$W_0=W_0^*=\ir,$$
$$W_k=A^k(W_0), $$
$$W_k^*=(A^*)^k(W_0). $$
Sequences $\{W_k\}_{k=0}^\infty$, $\{W^*_k\}_{k=0}^\infty$ are constructed by the Programmed Iteration Method for the initial and transformed problems respectively. Further let $\mathfrak{W}$ and $\mathfrak{W}^*$ denote the sets of approach problem solvability for the initial and transformed games respectively.

\begin{Th}
Let $M$ be a controllability set of the system (\ref{g_control}) and the target set $M^*=\{\vartheta\}\times F$. If for all $u\in P$, $v\in Q$, $\omega\in \Omega$, $\tau',\tau''\geq 0$ flows $\mathcal{F}^{\tau'}_{u,v}$ and $\mathcal{G}^{\tau''}_\omega$ commute: \begin{equation}\label{com_eq_one}
    \mathcal{F}^{\tau'}_{u,v}\circ \mathcal{G}^{\tau''}_\omega=\mathcal{G}^{\tau''}_\omega\circ\mathcal{F}^{\tau'}_{u,v},
\end{equation} then the following statements are fulfilled
\begin{enumerate}
\item\label{p_mpi} ${W}_k={W}^*_k $ for all $k\in \mathbb{N}\cup \{0\}$;
\item\label{p_eq} $\mathfrak{W}=\mathfrak{W}^*$;
\item\label{p_isaacs} if system (\ref{sys}) satisfies Isaacs condition, then the system
(\ref{eq_sys}) satisfies the Isaacs condition too.
\end{enumerate}
\end{Th}

The proof of this theorem is given in the end of section 4.

Note that if $f(\cdot,u,v)$ and $g(\cdot,\omega)$ are the smooth vector fields then the condition (\ref{com_eq_one}) may be written with the help of commutator of vector fields \cite{Agrachev}: $[\cdot,\cdot]$ $$[f(\cdot,u,v),g(\cdot,\omega)]=0\probel\forall u\in P\ \ \forall v\in Q\ \ \forall \omega\in \Omega. $$

Let us consider some examples. At first we consider the $M$-approach problem for system (\ref{sys}) when $M=[0,\vartheta]\times F$. In this case one can choose $\Omega=\{\omega\}$, $g(x,\omega)\equiv 0. $ Obviously, $M$ is controllability set of control system $g(x,\omega)$ and target set $M^*=\{\vartheta\}\times F$. Further, in this case $[f(x,u,v),g(x,\omega)]=0$. Thus the initial approach is equivalent to the $M^*$-approach for the conflict controlled system
$$\dot{x}=u_0\cdot f(x,u,v),\probel x\in \mathbb{R}^n,\probel u_0\in \{0,1\},\probel u\in P,\probel v\in Q. $$

The transformation in this case first was suggested in \cite{Mitchel}. In mentioned paper the statement 2 is proved for the case of the cylindrical target set. Actually, for differential games with simple motions and cylindrical target case the transformation were obtained by A.I.Subbotin \cite{Subb_book}. Moreover A.I. Subbotin got the formula of differential game solution in this case which is analogous to the formula obtained by B.N.~Pshenichny \cite{Pshen}.

Now we shall consider the problem of pointing of material point at sinking island with zero velocity.
Let $\vartheta=1$. Consider the conflict controlled system
$$ \left\{
\begin{array}{cc}
  \dot{y}= & z \\
  \dot{z}= & h(u,v).
\end{array}\right.
$$ Here $y$ and $z$ are 3d vectors, $u\in P$, $v\in Q$.
Suppose $$M=\{(t,y,z):t\in [0,1], \|y\|\leq 1-t,\probel z=0\}.$$ Hence $F$ contains only 6d zero. Choose $\Omega=\{\omega\in\mathbb{R}^3:\|\omega\|\leq 1\}$, $g(x,\omega)=g(\omega)=\omega$. It is easy to prove that $M$ is controllability set of chosen system $g(x,\omega)$, $\omega\in\Omega$. We have (see \cite{Agrachev})
$$[f(y,z,u,v),g(\omega)]=
\left(\begin{array}{cc}
  \mathbf{0} & \mathbf{0} \\
  \mathbf{0} & \mathbf{0}
\end{array}
\right)
\left(\begin{array}{c}
  z \\
  h(u,v)
\end{array}
\right)-\left(
\begin{array}{cc}
  \mathbf{0} & E \\
  \mathbf{0} & \mathbf{0}
\end{array}
\right)
\left(\begin{array}{c}
  g(\omega) \\
  0
\end{array}
\right)=0.
 $$ Here $\mathbf{0}$ means zero $3\times 3$ matrix, $E$ means identity $3\times 3$ matrix.

Therefore the problem of pointing of material point at the sinking island is equivalent to the problem of approach on the point $y=z=0$ at the moment $t=1$ for the system
$$ \left\{
\begin{array}{cc}
  \dot{y}= & \nu \cdot z+(1-\nu)\cdot g(\omega) \\
  \dot{z}= & \nu \cdot h(u,v).
\end{array}\right.
$$ In this system first player governs by the variables $\nu\in \{0,1\}$, $u\in P$, $\omega\in \Omega=\{y\in\mathbb{R}^3:\|y\|\leq 1\}$. The second player governs by the variable $v\in Q$.

\section{Some Properties of Program Absorption Operator}
Let us introduce one more property of sets. Let $E\subset \ir$. We say that $E$ {\it decreases} by section relative to control system $g(x,\omega)$, $\omega\in \Omega$, if for all $(t_*,x_*)$, $t\in [0,t_*]$ and $\sigma\in\mathcal{R}_\Omega$ the following inclusion holds: $\varphi_g(t,t_*,x_*,\sigma)\in E[t]$.

\begin{Lm}\label{lm_keep}
Let $E\subset\ir$. If $M\subset E$ and $E$ decreases by section relative to control system $g(x,\omega)$, $\omega\in\Omega$, then $A^*(E)$ posses these properties.
\end{Lm}
\begin{proof}
At first we prove that $A^*(E)$ decreases by section relative to control system $g(x,\omega)$, $\omega\in\Omega$. Let $(t_*,x_*)\in A^*(E)$, $t\in [0,t_*]$, $\sigma\in \mathcal{R}_\Omega$. Our purpose is to prove that $\varphi_g(t,t_*,x_*,\sigma)\in (A^*(E))[t]$. Since $(t_*,x_*)\in A^*(E)\subset E$, we have \begin{equation}\label{tau_e}
(\tau,\varphi_g(\tau,t_*,x_*,\sigma))\in E\probel\forall \tau\in [t,t_*].
\end{equation}
For all $v\in Q$ there exists measure $\mu\in\mathcal{R}_{P^*}$ such that  $\phi^*(\xi,t_*,x_*,\mu,v)\in M[\xi]$ for some $\xi\in [t_*,\vartheta]$ and for all  $\tau\in [t_*,\xi]$ the inclusion $\phi^*(\tau,t_*,x_*,\mu,v)\in E[\tau]$ is fulfilled.

Also there exists measure $\tilde{\sigma}\in\mathcal{R}_{P^*}$ such that
$$\int_{[0,\vartheta]\times P^*}\psi(t,\omega)\tilde{\sigma}(d(t,\nu,\omega,u))=\int_{[0,\vartheta]\times\Omega}\psi(t,\omega)\sigma(d(t,\omega)) $$ for all $\psi\in C([0,\vartheta]\times\Omega)$.

Let $\tilde{\mu}$ be a measure such that $$\int_{[t,\vartheta]\times P^*}\psi(t,u^*)\tilde{\mu}(d(t,u^*))=\int_{[t,t_*]\times P^*}\psi(t,u^*)\tilde{\sigma}(d(t,u^*))+\int_{[t_*,\vartheta]\times P^*}\psi(t,u^*){\mu}(d(t,u^*)). $$
Put $\bar{x}=\varphi_g(t,t_*,x_*,\sigma)$. Note that
\begin{equation}\label{varphi_chain}
\varphi_g(\tau,t_*,x_*,\sigma)=\varphi_g(\tau,t,\bar{x},\sigma)=\phi^*(\tau,t,\bar{x},\tilde{\sigma},v)= \phi^*(\tau,t,\bar{x},\tilde{\mu},v).
\end{equation}
Combining (\ref{tau_e}) and (\ref{varphi_chain}) we get
\begin{equation}\label{tau_1_incl}
\phi^*(\tau,t,\bar{x},\tilde{\mu},v)\in E[\tau].
\end{equation}
Also we have $x_*=\varphi_g(t_*,t,\bar{x},\sigma)=\phi^*(\tau,t,\bar{x},\tilde{\mu},v).$ Therefore $$\phi^*(\tau,t_*,x_*,\mu,v)=\phi^*(\tau,t_*,x_*,\tilde{\mu},v)=\phi^*(\tau,t,\bar{x},\tilde{\mu},v).$$ Since $(t_*,x_*)\in A^*(E)$ we claim that $\phi^*(\xi,t,\bar{x},\tilde{\mu},v)\in M[\xi]$ and for all $\tau\in [t_*,\xi]$ $\phi^*(\tau,t,\bar{x},\tilde{\mu},v)\in E[\tau]$. Combining this with inclusion (\ref{tau_1_incl}) we obtain $(t,\bar{x})\in A^*(E)$.

The inclusion  $M\subset A^*(E)$ follows from the definition of $A^*$, representation (\ref{M_control}), and formula (\ref{f_star}).
\end{proof}

Now let $u^*(\cdot):[0,\vartheta]\rightarrow P^*$ be a piecewise function. Suppose $t_*,\xi\in [0,\vartheta]$, $t_*\leq \xi$. The half-interval $[t_*,\xi)$ can be represented as union of half-intervals $[\xi_{i-1},\xi_i)$, $i=\overline{1,k}$, such that $u^*(\theta)=u_i$, $\theta\in [\xi_{i-1},\xi_i)$. Here $u_1,\ldots, u_k$ are the elements of control space $P^*$; $u^*_i=(\nu_i,u_i,\omega_i)$, $\nu_i\in \{0,1\}$, $u_i\in P$, $\omega_i\in \Omega$. Put $$J'\triangleq\{i:\nu_i=1\}=\{r_1,\ldots, r_l\},$$
$$J''\triangleq\{i:\nu_i=0\}. $$

Denote $\hat{\xi}_0\triangleq t_*$, $\hat{\xi}_i\triangleq \hat{\xi}_j+\tau_{r_j}$,
 \begin{equation}\label{xi_def}
 \bar{\xi}\triangleq \xi_{r_l}.
\end{equation}

Define the piecewise control $u(\cdot):[t_*,\bar{\xi})\rightarrow P$ by the rule
\begin{equation}\label{u_def}
u(t)=u_{r_j}, \probel\mbox{for } t\in [\hat{\xi}_{r_{j-1}},\hat{\xi}_{r_j}).
\end{equation}

Suppose that $t\in [t_*,\bar{\xi}]$. Either there exists $j$ such that $t\in [\hat{\xi}_{r_{j-1}},\hat{\xi}_{r_{j}})$, or $t=\bar{\xi}$. In the first case put $\gamma(t)\triangleq \xi_{r_j-1}+t-\hat{\xi}_{r_{j-1}}$, in the second case put $\gamma(t)\triangleq \bar{\xi}$.
Moreover denote $J''_{t}=\{i\in J'':\xi_i<\gamma(t)\}=\{s_1,\ldots,s_m\}$.
\begin{Lm}\label{lm_equality} There exists piecewise control $\omega(\cdot):[0,\vartheta]\rightarrow \Omega$ such that
$$\mathbf{x}^*(\gamma(t),t_*,x_*,u^*(\cdot),v)=x_{g}(\gamma(t),t,\mathbf{x}(t,t_*,x_*,u(\cdot),v),\omega(\cdot))$$ for all $t\in [t,\bar{\xi}]$.
\end{Lm}
\begin{proof}
The following representation is fulfilled: $$\mathbf{x}^*(\gamma(t),t_*,x_*,u^*(\cdot),v)=\mathcal{F}^{*\probel \gamma(t)-\xi_{r_j-1}}_{u_{s_j}^*,v}\circ\ldots\circ\mathcal{F}^{*\probel\tau_i}_{u_i^*,v}\circ \ldots\circ \mathcal{F}^{*\probel \tau_1}_{u_1^*,v}(x_*). $$

We have $\mathcal{F}^{*\probel \gamma(t)-\xi_{r_j-1}}_{u_{r_j}^*,v}=\mathcal{F}^{ t-\hat{\xi}_{r_{j-1}}}_{u_{r_j},v}$. Also for $i\in J'$, $\mathcal{F}^{*\probel\tau_i}_{u_i^*,v}=\mathcal{F}^{\tau_i}_{u_i,v}$, for $i\in J''$ $\mathcal{F}^{*\probel\tau_i}_{u_i^*,v}=\mathcal{G}^{\tau_i}_{\omega_i}$. The flows $\mathcal{F}^{\tau'}_{u,v}$ and $\mathcal{G}^{\tau''}_{\omega}$ commute by the assumption of Theorem. Therefore $$\mathbf{x}^*(\gamma(t),t_*,x_*,u^*(\cdot),v)=\mathcal{G}^{\tau_{s_m}}_{\omega_m}\circ\ldots\circ \mathcal{G}^{\tau_{s_1}}_{\omega_1}\circ\mathcal{F}^{ t-\hat{\xi}_{r_{j-1}}}_{u_{r_j}^*,v}\circ\ldots\circ  \mathcal{F}^{ \tau_{r_1}}_{u_{r_1},v}(x_*). $$ This completes the proof.
\end{proof}

\begin{Lm}\label{lm_eq_a}
 Let $E$ decrease by sections relative to $g(x,\omega)$, $\omega\in \Omega$. Also let $M\subset E$. Then $A(E)=A^*(E)$.
\end{Lm}
\begin{proof}
At first we shall prove the inclusion $A(E)\subset A^*(E)$. Let $(t_*,x_*)\in A(E)$. By definition of operator $A$ we have that for all  $v\in Q$ there exist measure $\mu\in\mathcal{R}_P$ and moment $\xi\in [t_*,\vartheta]$ such that $\phi(\xi,t_*,x_*,\mu,v)\in M[\xi]$ and for all   $t\in [t_*,\xi]$  the  inclusion  $\phi (t,t_*,x_*,\mu,v)\in E[t]$ holds. Since the set $M$ is a controllability set of control system $g(x,\omega)$, $\omega\in\Omega$, there exists measure $\sigma\in \mathcal{R}_\Omega$ such that $\varphi_g(t,\xi,\phi(\xi,t_*,x_*,\mu,v),\sigma,v)\in M[t]$, $t\in [\xi,\vartheta]$. By the Riss theorem one can choose measures $\hat{\mu}\in \mathcal{R}_{P^*}$ � $\tilde{\sigma}\in \mathcal{R}_{P^*}$ such that $$\int_{[0,\vartheta]\times P}\psi(t,u)\mu(d(t,u))=\int_{[0,\vartheta]\times P^*}\psi(t,u)\hat{\mu}(d(t,\nu,u,\omega)) $$ for all $\psi\in C([0,\vartheta]\times P)$, and
$$\int_{[0,\vartheta]\times \Omega}\psi(t,\omega)\sigma(d(t,\omega))=\int_{[0,\vartheta]\times P^*}\psi(t,\omega)\tilde{\sigma}(d(t,\nu,u,\omega)) $$ for all $\psi\in C([0,\vartheta]\times\Omega)$.

Let $\beta\in\mathcal{R}_{P^*}$ be a measure such that for all functions $\psi\in C([0,\vartheta]\times P^*)$  the following equalities hold:
$$\int_{[0,\xi]\times P^*}\psi(t,\nu,u,\omega)\beta(d(t,\nu,u,\omega))=\int_{[0,\xi]\times P^*}\psi(t,\nu,u,\omega)\hat{\mu}(d(t,\nu,u,\omega)), $$
$$\int_{[\xi,\vartheta]\times P^*}\psi(t,\nu,u,\omega)\beta(d(t,\nu,u,\omega))=\int_{[\xi,\vartheta]\times P^*}\psi(t,\nu,u,\omega)\tilde{\sigma}(d(t,\nu,u,\omega)). $$

We have $$\phi^*(t,t_*,x_*,\beta,v)=\left\{
\begin{array}{cc}
  \phi(t,t_*,x_*,\mu,v),\sigma,v), & t\in [t_*,\xi]. \\
  \varphi_g(t,\xi,\phi(\xi,t_*,x_*,\mu,v),\sigma,v), & t\in [\xi,\vartheta].
\end{array}
\right. $$
Since $M[t]\subset E[t]$ and $M[\vartheta]=F$ we claim that $(\vartheta,\phi^*(\vartheta,t_*,x_*,\beta,v)\in \{\vartheta\}\times F$ and for all $t\in [t_*,\vartheta]$ the inclusion $\phi^*(t,t_*,x_*,\beta,v)\in E[t]$ holds.

Now we shall prove the inclusion $A^*(E)\subset A(E)$.

Choose $(t_*,x_*)\in A^*(E)$. Let $C>0$ be a number such that for all $t_1,t_2\in [0,\vartheta]$, $t_2\leq t_1$, $x',x''\in \mathcal{G}$, $\sigma\in\mathcal{R}_\Omega$ the inequality $$\|\varphi_g(t_2,t_1,x',\sigma)-\varphi_g(t_2,t_1,x'',\sigma)\|\leq C\|x'-x''\|$$ is fulfilled. Here $\mathcal{G}$ is the accessibility set  from $[0,\vartheta]\times \{x_*\}$ under the action of control system $\dot{x}=f^*(x,\nu,u,\omega,v),\probel \nu\in \{0,1\}\probel u\in P,\probel \omega\in\Omega,\probel v\in Q$.

Inclusion $(t_*,x_*)\in A^*(E)$ means that for all $v\in Q$ one can choose a measure $\beta\in \mathcal{R}_{P^*}$ such that $\phi^*(\vartheta,t_*,x_*,\beta,v)\in F$ and for all $t\in [t_*,\vartheta]$ the inclusion $\phi(t,t_*,x_*,\beta,v)\in E[t]$ is fulfilled.
Further there exists a sequence of piecewise controls $\{\zeta^\alpha(\cdot)\}_{\alpha=1}^\infty$, $\zeta^\alpha(\cdot):[t_*,\vartheta]\rightarrow P^*$ such that $$\varepsilon^\alpha\triangleq\sup_{t\in [t^*,\vartheta]}\|\mathbf{x}^*(t,t_*,x_*,\zeta^\alpha(\cdot),v)-\phi^*(t,t_*,x_*,\beta,v)\|\rightarrow 0, \alpha\rightarrow\infty.$$ Let us consider the sequence of controls $\{u^\alpha\}_{\alpha=1}^\infty$ and sequence of moments $\{\xi^\alpha\}_{\alpha=1}^\infty$, for those the elements $\xi^\alpha$ and $u^\alpha$ are defined by the rules (\ref{xi_def}) and (\ref{u_def}) respectively. Further for each $\alpha$ the function $\gamma^\alpha(\cdot)$ is well defined.

There exists a subsequence $\{\alpha_k\}$ such that $\xi^{\alpha_k}\rightarrow\xi$, $\mu_{\zeta^{\alpha_k}}\rightharpoondown \mu$. Without loss of generality it can be assumed that subsequence $\{\alpha_k\}$ coincides with sequence $\{\alpha\}$.

We have $\gamma^\alpha(\xi^\alpha)=\vartheta$ for all $\alpha\in\mathbb{N}$. Using lemma \ref{lm_equality} we obtain that   for some control $\omega^\alpha(\cdot)$ the following equality holds: $$\mathbf{x}^*(\vartheta,t_*,x_*,\zeta^\alpha(\cdot),v)= x_g(\vartheta,\xi,\mathbf{x}(\xi,t_*,x_*,u^\alpha(\cdot),v),\omega^\alpha(\cdot)).$$ This is equivalent to the equality $$\mathbf{x}(\xi,t_*,x_*,u^\alpha(\cdot),v)= x_g(\xi,\vartheta,\mathbf{x}^*(\vartheta,t_*,x_*,\zeta^\alpha(\cdot),v),\omega^\alpha(\cdot)). $$ Therefore \begin{multline*}
\|\mathbf{x}(\xi^\alpha,t_*,x_*,u^\alpha(\cdot),v)-x_g(\xi,\vartheta,\phi^*(\vartheta,t_*,x_*,\beta,v),\omega^\alpha(\cdot))\|\leq \\ \leq \|x_g(\xi^\alpha,\vartheta,\mathbf{x}^*(\vartheta,t_*,x_*,\zeta^\alpha(\cdot),v),\omega^\alpha(\cdot))-x_g(\xi,\vartheta,\phi^*(\vartheta,t_*,x_*,\beta,v),\omega^\alpha(\cdot))\| \leq C\varepsilon^\alpha.
\end{multline*} Since $M$ is controllability set of control system $g(x,\omega)$, $\omega\in\Omega$, and target set $\{\vartheta\}\times F$ we have $$(\xi^\alpha,x_g(\xi^\alpha,\vartheta,\phi^*(\vartheta,t_*,x_*,\beta,v),\omega^\alpha(\cdot)))\in M.$$ Consequently  $$(\xi,\mathbf{x}(\xi,t_*,x_*,\mu,v))\in M.$$

Now let $t\in [t_*,\xi]$. For sufficient large $\alpha$ the inequality $t\leq\xi^\alpha$ is fulfilled. We have $\phi^*(\gamma^\alpha(t),t_*,x_*,\beta,v)\in E[\gamma(t)]$. Using lemma (\ref{lm_equality}) we get, that for some control $\omega(\cdot)$ the following equality holds: $$\mathbf{x}(t,t_*,x_*,u^\alpha(\cdot),v)= x_g(t,\gamma^\alpha(t),\mathbf{x}^*(\gamma^\alpha(t),t_*,x_*,\zeta^\alpha(\cdot),v),\omega(\cdot)). $$

Therefore \begin{multline*}
    \|\mathbf{x}(t,t_*,x_*,u^\alpha(\cdot),v)-x_g(t,\gamma^\alpha(t),
\phi^*(\gamma^\alpha(t),t_*,x_*,\beta,v),\omega(\cdot))\|\leq\\
\leq\|x_g(t,\gamma^\alpha(t),\mathbf{x}^*(\gamma^\alpha(t),t_*,x_*,\zeta^\alpha(\cdot),v),\omega(\cdot))-x_g(t,\gamma^\alpha(t),
\phi^*(\gamma^\alpha(t),t_*,x_*,\beta,v),\omega(\cdot))\|\leq C\varepsilon^\alpha.
\end{multline*}
Since $E$ decreases by sections relative to control system $g(x,\omega)$, $\omega\in \Omega$ and the inclusion $\phi^*(\gamma^\alpha(t),t_*,x_*,\beta,v)\in E[\gamma^\alpha(t)]$ is fulfilled, we conclude that $x_g(t,\gamma^\alpha(t),
\phi^*(\gamma^\alpha(t),t_*,x_*,\beta,v),\omega(\cdot))\in E[t]$.

Since the choose of $v\in Q$ is arbitrary, we obtain $(t_*,x_*)\in A(E)$. Thus it is established that $A^*(E)\subset A(E)$.

\end{proof}

\begin{proof}[The proof of Main Theorem]
The statement 1 follows from lemmas 1 and 3 since $M\subset \ir$ and the set $W_0=W_0^*=\ir$ decreases by sections relative to control system $g(x,\omega)$, $\omega\in \Omega$.

The statement 2 obviously follows from the statement 1 and representation of the solvability set (\ref{Chentsov_formula}).

Now we shall prove the statement 3. The inequality
$$\max_{v\in Q}\min_{(\nu,u,\omega)\in P^*}\langle s, f^*(x,\nu,u,\omega,v)\rangle\leq \min_{(\nu,u,\omega)\in P^*}\max_{v\in Q}\langle s, f^*(x,\nu,u,\omega,v)\rangle $$ is obvious. Thus, we need to prove opposite inequality.

Denote  a saddle point in the small game for system (\ref{sys}) by $(u_*,v_*)$. Let us consider two cases. At first suppose that $$
\min_{(\nu,u,\omega)\in P^*}\max_{v\in Q}\langle s, f^*(x,\nu,u,\omega,v)\rangle=\min_{u\in P,\omega\in \Omega}\max_{v\in Q}\langle s, f^*(x,1,u,\omega,v)\rangle.
$$ In particular it means that  $$\max_{v\in Q}\min_{u\in P}\langle s, f(x,u,v)\rangle=\min_{u\in P}\max_{v\in Q}\langle s, f(x,u,v)\rangle\leq \min_{\omega\in\Omega}\langle s,g(x,\omega)\rangle. $$ In this case
\begin{multline*}\min_{(\nu,u,\omega)\in P^*}\max_{v\in Q}\langle s, f^*(x,\nu,u,\omega,v)\rangle=\min_{u\in P}\max_{v\in Q}\langle s, f(x,u,v)\rangle=\max_{v\in Q}\min_{u\in P}\langle s, f(x,u,v)\rangle=\\=\min_{u\in P}\langle s, f(x,u,v_*)=\min\{\min_{u\in P}\langle s, f(x,u,v_*),\min_{\omega\in\Omega}\langle s,g(x,\omega)\rangle\} = \\ =
\min_{(\nu,u,\omega)\in P^*}\langle s, f^*(x,\nu,u,\omega,v_*)\rangle\leq \max_{v\in Q}\min_{(\nu,u,\omega)\in P^*}\langle s, f^*(x,\nu,u,\omega,v)\rangle.
\end{multline*}

Now suppose that $$
\min_{(\nu,u,\omega)\in P^*}\max_{v\in Q}\langle s, f^*(x,\nu,u,\omega,v)\rangle=\min_{u\in P,\omega\in \Omega}\max_{v\in Q}\langle s, f^*(x,0,u,\omega,v)\rangle=\min_{\omega\in\Omega}\langle s,g(x,\omega)\rangle.
$$ Then \begin{multline*}
\min_{\omega\in\Omega}\langle s,g(x,\omega)\rangle= \min\{\min_{\omega\in\Omega}\langle s,g(x,\omega)\rangle,\min_{u\in P}\max_{v\in Q}\langle s,f(x,u,v)\rangle\}=\\=\min\{\min_{\omega\in\Omega}\langle s,g(x,\omega)\rangle,\min_{u\in P}\langle s,f(x,u,v_*)\rangle\}=\min_{(\nu,u,\omega)\in P^*}\langle s,f(x,\nu,u,\omega,v_*)\rangle\leq\\\leq \max_{v\in Q}\min_{(\nu,u,\omega)\in P^*}\langle s,f(x,\nu,u,\omega,v)\rangle.
\end{multline*}

\end{proof}

\noindent{\bf Acknowledgments.} The author is grateful to member-correspondent of RAS A.G. Chentsov for the great attention to this paper.


\begin{thebibliography}{99}
\bibitem{isaacs}Isaacs R., {\it Differential Games}, John Wiley and Sons, 1965.

\bibitem{NN_PDI_en} Krasovskii N.N., Subbotin A.I., {\it Game-theoretical control problems}, Springer-Verlag, 1988, 517 p.
\bibitem{Subb_Ch}  Subbotin N.N., Chentsov A.G. {\it Optimization
of a Guarantee in Control Problem,} Nauka, Moscow, 1980 (in
Russian).

\bibitem{Ch_dan_75} Chentsov A.G., On the structure of a game problem of convergence, {\it  Soviet Math. Dokl.}, {\bf 16} (1975), 1404--1406.
\bibitem{Math_sb} Chentsov A.G., On a game problem of converging at a given instant time,
{\it Math. USSR Sbornic,} {\bf 28}, 3 (1976), 353--376.

\bibitem{Chentsov_76} Chentsov A.G., On a game problem of guidance
{\it  Soviet Math. Dokl.}, {\bf 17}, 1 (1976), 73--77.

\bibitem{Subb_book} Subbotin A.I, {\it Generalized Solutions of First-Order PDEs. The Dynamical
Optimization Perspective.} Birkhauser, Boston, 1995.

\bibitem{Mitchel} Mitchel I.M., Bayen A.M., Tomlin C.J., A Time-Depend Hamilton-Jacobi Formulation of Reachable Sets for
Continuous Dynamic Games // {\it IEEE Transaction on Automatic Control}, {\bf 50},  7 (2005), 
947--957.

\bibitem{Agrachev} Agrachev A.A. Sachkov Yu.L., {\it Control theory from the geometric viewpoint}. Springer-Verlag, Berlin, Heidelberg, New York, Tokyo, 2004. 

\bibitem{Chen_Dep} Chentsov A.G., {\it The programmed iteration method for a differential
pursuit-evasion game}, Dep. in VINITI, 1933-79, Sverdlovsk, 1979 (in
Rusian).

\bibitem{NN_Cyber} Krasovskii N.N., Differential Game of approach-deviation - I, {\it Izv. AN USSR (Tech.
cybernetics)}, {\bf 2}  (1973), 3--18, (in Rusian).

\bibitem{Pshen} {\it Pschenitchny B.N.} $\varepsilon$-strategies in differential games, in: Topics in Differential Games, North Holland, New York-London-Amsterdam (1973), pp. 45-99.
\end{thebibliography}
\end{document}